\documentclass[11pt]{article}
\usepackage{amsfonts,amssymb,graphicx}
\usepackage{amsfonts}
\usepackage{amsmath}
\usepackage{amstext}
\usepackage{amsthm}
\usepackage{xspace}

\newtheorem{thm}{Theorem}[section]
\newtheorem{cor}[thm]{Corollary}
\newtheorem{lem}[thm]{Lemma}
\newtheorem{prop}[thm]{Proposition}

\theoremstyle{definition}

\numberwithin{equation}{section}

\newcommand{\set}[1]{\left\{#1\right\}}
\newcommand{\R}{\mathbb R}

\newcommand{\C}{{\bf \mathcal{C}}}
\newcommand{\res}{{\bf \mathcal{S}}}
\newcommand{\B}{{\bf \mathcal{B}}}
\newcommand{\Z}{\mathbb Z}
\newcommand{\Q}{\mathbb Q}
\newcommand{\N}{\mathbb N}

\newcommand{\virg}[1]{``#1"}

\begin{document}

\title{A renormalization approach to irrational rotations}

\author{Claudio Bonanno\footnote{Dipartimento di Matematica Applicata,
Universit\`a di Pisa, via F. Buonarroti 1/c, 56127 Pisa (Italy),
email: $<$bonanno@mail.dm.unipi.it$>$;} \and Stefano
Isola\footnote{Dipartimento di Matematica e Informatica,
Universit\`a di Camerino, via Madonna delle Carceri, 62032
Camerino (Italy), e-mail: $<$stefano.isola@unicam.it$>$}}

\maketitle
% ----------------------------------------------------------------

\begin{abstract}

\noindent We introduce a renormalization procedure which allows us
to study in a unified and concise way different properties of the
irrational rotations on the unit circle $\beta \mapsto
\set{\alpha+\beta}$, $\alpha \in \R\setminus \Q$. In particular we
obtain sharp results for the diffusion of the walk on $\Z$
generated by the location of points of the sequence $\{n\alpha
+\beta\}$ on a binary partition of the unit interval. Finally we
give some applications of our method.

\end{abstract}
% ----------------------------------------------------------------

2000 {\sl Mathematics Subject Classification}: 37E05, 37A25, 37A45

{\sl Key words and phrases}: Irrational rotations, continued
fractions, renormalization, diffusion.

\section{Introduction} \label{sec:intro}

Irrational rotations on the unit circle $S^1\cong [0,1]/(0=1)$ are
isometric transformations defined by
$$[0,1) \ni \beta \mapsto \set{\alpha + \beta} \in [0,1)$$
where $\alpha \in \R \setminus \Q$ is the angle of rotation and
$\set{\cdot}$ denotes the fractional part of a real number. It is
well known that the Lebesgue measure on the unit interval is the
unique (and thus ergodic) invariant probability measure for these
transformations. However, some ergodic properties of the
rotations, such as recurrence rates, waiting times and some limit
laws, are known to depend on the arithmetic properties of the
rotation angle $\alpha$ (see \cite{bgi}, \cite{chaz},
\cite{coelho}, \cite{kim}, \cite{ks}). These, in turn, are encoded
in its continued fraction expansion and it is starting from this
expansion that we introduce, borrowing it from dynamical systems,
a {\sl renormalization procedure} which allows us to study several
relevant properties of the orbit $(\set{n\alpha + \beta})$. The
possibility of studying all these properties through the same
approach is one of the main motivations of this paper. Indeed, the
paper is self-contained and we obtain sharp results using
efficient and concise arguments based solely on the
renormalization procedure.

A renormalization approach to circle maps was introduced in
\cite{lanford}, and used for example in \cite{coelho} to estimate
limit laws of entrance times for irrational rotations. We remark
that the renormalization procedure we use is quite different from
the one in \cite{lanford}, both in the construction and in the
spirit. Moreover, we think that our approach can be used
successfully to give sharp estimates on some more ergodic
properties and limit laws than those considered in this paper (see
\cite{bon}).

\noindent In this paper, we study the distribution properties of
the sequences $x_n:=n\alpha + \beta$, $n=0,1,2,\dots$. We recall
that a sequence $(y_n)$ of real numbers is said to be {\sl
uniformly distributed modulo $1$}, if for any real number
$0<\gamma \leq 1$ we have
$$
\lim_{n\to \infty} \frac 1 n \ \sum_{r=0}^{n-1} \
\chi_{[0,\gamma)} (\set{y_r})=\gamma
$$
where $\chi_{A}$ denotes the indicator function of the set $A$.

An interesting characterisation of the distribution properties of
a sequence $(y_n)$ can be obtained as follows: take the partition
of the unit interval given by $\set{[0,\frac 1 2),[\frac 1 2,1)}$
and construct the walk on $\Z$ which starts at the origin at time
$0$ and at time $n+1$ moves one step to the right if $\set{y_n}
\in [0,\frac 1 2)$, one step to the left otherwise.  After $n+1$
steps, the position of the walker is given by

\begin{equation} \label{def:somme}
S_n=\sum_{r=0}^{n}s(\set{y_r})
\end{equation}

$$s(y):=2\, \chi_{[0,\frac 1 2)}(\set{y})-1.$$

\noindent Clearly, if $(y_n)$ is uniformly distributed modulo $1$
then $|S_n|=o(n)$ as $n\to \infty$ and the better is the uniform
distribution of $(y_n)$ the slower is the diffusion of the walk.
In particular, for the ideal distribution for which $\frac 1 n \
\sum_{r=0}^{n-1} \ \chi_{[0,\gamma)} (\set{y_r})=\gamma$ for all
$n$ we have $S_n=0$ for all $n$.

In \cite{isola} the growth of the quantity $S_n$ has been studied
for the sequence $(x_n)=(n\alpha + \beta)$ in the $L^\infty$ and
$L^2$ norms. In this paper we give an explicit formula for $S_n$
and obtain as a corollary growth estimates for $\beta=0$, and a
sharp result for the $L^\infty$ norm (see Corollary
\ref{cor:l-inf}), improving results in \cite{isola}.

\noindent For the sequence $x_n:= n\alpha + \beta$, with $\alpha$
irrational and $\beta\in (0,1)$, the indicator $S_n$ will be
denoted by $S_n(\alpha,\beta)$ ($S_n(\alpha)$ if $\beta=0$) to
stress its dependence on the arithmetical properties of the number
$\alpha$.

\noindent The paper is organized as follows. Notations and the
main ideas are settled in Section \ref{sec:rot} which includes a
preliminary analysis of several quantities that are needed for the
renormalization procedure. In particular, we study how the
different points of the sequence $(n\alpha)$ are organized
according to their integer part. This depends on $a_1$, the first
partial quotient in the continued fraction expansion of $\alpha$,
and in particular by its parity. The main result of this section
is Theorem \ref{cor:rkj}.

\noindent In Section \ref{sec:diffusione} we give an iterative
method to obtain an explicit expression for $S_n(\alpha)$ only in
terms of the coefficients of the continued fraction of $\alpha$.
This method is then used to obtain sharp estimates for
$S_n(\alpha)$ that have significantly different expressions
according to whether $a_1$ is even or odd. The main results are
the algorithm described in Proposition \ref{prop:massimi} and its
consequences described in Theorem \ref{cor:maxpari} as well as in
the subsequent examples. We point out that we obtain results also
for the minimum values of $S_n(\alpha)$, with the aim of giving
hints on the returns to zero (this point is analysed in
\cite{bon}).

\noindent Next, we consider the general case of the sums
$S_n(\alpha,\beta)$. Again we obtain an explicit expression for
these sums (Theorem \ref{prop:gen-diff}) in terms of the
coefficient of the continued fraction of $\alpha$ and of the
coefficient of the expansion of $\beta$ introduced in Proposition
\ref{prop:espansione}. In particular we obtain as a corollary a
sharp estimate on the $L^\infty$ norm of $S_n(\alpha,\beta)$ (see
Corollary \ref{cor:l-inf}).

\noindent Finally, we give some applications of our approach to
the Birkhoff Ergodic Theorem and to the \emph{discrepancy} (see
(\ref{def:disc})) of the sequence $(n\alpha)$. It is surprising
that these results easily follow from our renormalization
approach.

\section{Continued fractions and return sequences} \label{sec:rot}

\noindent For a given  number $\alpha \in (0,1)$ let us consider
its expansion in continued fraction \cite{pytheasfogg}
\begin{equation} \label{eq:cont-frac}
\alpha = {1\over \displaystyle a_1 +
{1\over \displaystyle a_2 + {1\over\displaystyle a_3 +\cdots }}}
\end{equation}
which we denote by $\alpha=  [a_1,a_2,a_3,\dots]$. The partial
quotients $a_h$ are positive integers and the expansion terminates
if and only if $\alpha$ is rational. If $\alpha$ is irrational its
\virg{fast} rational approximants are the numbers
\begin{equation}
\frac{p_n}{q_n} := [a_1,\dots, a_n]
\end{equation} which can be also be defined recursively by
\begin{equation} \label{eq:rel-approx}
{p_0\over q_0}={0\over 1},\quad {p_1\over q_1}={1\over a_1}
\quad\hbox{and}\quad {p_{n+1}\over q_{n+1}}={a_{n+1}\, p_n +
p_{n-1}\over a_{n+1}\, q_n + q_{n-1}},\quad n\geq 1.
\end{equation}
In the following we shall consider also the positive numbers
\begin{equation} \label{effenne}
f_n:=(-1)^n (q_n \alpha - p_n),\qquad n\ge 0.
\end{equation}
To a given $\alpha \in (0,1)$ we associate the rotation
$T_\alpha:X\to X$ of the unit circle $X=[0,1]/(0=1)$ given by
\begin{equation}
T_\alpha(\beta) := \set{\alpha+\beta}.
\end{equation}
One easily checks that the numbers $f_n$ determine the successive
closest returns of the orbit of a point $x$ to the point itself
(thus forming a monotonically decreasing sequence). Indeed, for
all $n\ge 1$ and for all $\beta \in X$, it holds
\begin{equation}
q_n = \min \set{ r> q_{n-1} \ :\ |x-T_\alpha^r(\beta)| <
|\beta-T_\alpha^{q_{n-1}}(\beta)|\,}
\end{equation}
and
\begin{equation}
f_n = |\beta-T_\alpha^{q_n}(\beta)|.
\end{equation}
The first numbers $f_n$ are given by $f_0=\alpha$, $f_1=1 -a_1
\alpha$, $f_2 = f_0 - a_2 f_1$ and more generally they satisfy the
recursion
\begin{equation} \label{recu}
f_{n+1}=f_{n-1} - a_{n+1} f_n,\qquad n\ge 1.
\end{equation}
Conversely, once the $f_n$ are known the partial quotients can be
obtained as:
\begin{equation} \label{recu1}
a_{n+1} = \max \set{h \ge 1 \ :\ h f_n < f_{n-1}},\qquad n\geq 0,
\end{equation}
with the position $f_{-1}=1$. This yields in particular
$$a_1 = \max \set{h \ge 1 \ :\ h \alpha<1}$$
$$a_2 = \max \set{h \ge 1 \ :\ h(1-a_1 \alpha) < \alpha}.$$

\vskip 0.2cm
\noindent Note that $\lfloor r\alpha \rfloor=0$ for
all $0\le r \le a_1$ and $\lfloor (a_1+1)\alpha \rfloor=1$. In the
sequel we shall study the behaviour of the sum $S_n(\alpha)$ by
looking at the values of $s(r\alpha)$ with $\lfloor r\alpha
\rfloor$ constant. To this end we introduce the following
quantities. Set
\begin{equation} \label{eq:rk}
r_k := \min \set{r\ge 0 \ :\ \lfloor r\alpha \rfloor=k},\quad k\geq 0.
\end{equation}
In terms of $T_\alpha$ this is the least number of iterates of $0$
needed to make $k$ ``turns'' of the circle $X$. Set moreover
\begin{equation} \label{eq:tk}
t_k:= \# \set{r\ge 0 \ :\ \lfloor r\alpha \rfloor =k}, \quad k\geq 0,
\end{equation}
which is the number of  $T_\alpha$-iterates of $0$ which are all
lying ``within the same circle'' after having turned the circle
$k$ times.

\noindent One sees that $t_0 = a_1 +1$, and it is not difficult to
realize that for all $k\ge 1$, $t_k$ is equal to either $a_1$ or
$a_1+1$. More precisely we have

\begin{lem} \label{lemmino}
$t_k$ is either equal to $a_1 +1$ or $a_1$, according to
whether $\set{r_k \alpha}$ is smaller or bigger than $f_1$, respectively.
\end{lem}

\noindent {\sl Proof.} Let $\lfloor r_k \alpha \rfloor=k$, then
$0\le \{ r_k \alpha\} \le \alpha$ and $\{ r_k \alpha\} + (a_1 -1)
\alpha \le a_1 \alpha < 1$, hence $t_k \ge a_1$. Moreover $\{ r_k
\alpha\} + (a_1 +1) \alpha \ge (a_1+1)\alpha >1$, hence $t_k \le
a_1+1$. Finally $t_k = a_1+1$ if and only if $\{ r_k \alpha\} +
a_1 \alpha <1$, that is if and only if $\{ r_k \alpha\} < (1- a_1
\alpha)=f_1$. $\Box$

\vskip 0.5cm
\noindent Starting our analysis from $r=0$, we notice
that $r_1=a_1+1$ and $\set{r_1 \alpha}=\alpha-f_1=f_0-f_1$. Then
$t_1 = a_1+1$ if and only if $f_0-f_1 < f_1$, that is if and only
if $f_2=f_0-f_1$, that is $a_2=1$ (cfr. (\ref{recu1})). If instead
$a_2>1$ then $\set{r_2 \alpha}= \set{r_1 \alpha}-f_1 = f_0-2f_1$
and proceeding recursively  $\set{r_k \alpha}=f_0 -kf_1 > f_1$ for
all $1\le k < a_2$, hence $t_k = a_1$ for all $1\le k < a_2$. On
the other hand $\set{r_{a_2} \alpha}=f_0 - a_2 f_1 = f_2<f_1$,
whence $t_{a_2}=a_1+1$ and $r_{a_2}=q_2$, the denominator of the
second \virg{fast} rational approximant of $\alpha$.

\noindent Let us denote by $(r_{k_j})$ the sub-sequence of $(r_k)$
such that $t_{k_j} = a_1+1$ for all $j\ge 0$. So far we have
showed that $k_0=0$, $r_{k_0}=0$ and $k_1=a_2$, $r_{k_1}=q_2$. We
now investigate the following terms of $(r_{k_j})$. Let $(g_j)$
denote the sequence of \virg{gaps} between subsequent elements of
$(r_{k_j})$:
\begin{equation} \label{eq:gj}
g_{j} := r_{k_j} - r_{k_{j-1}},\qquad j\geq 1.
\end{equation}
Given the irrational number $\alpha=[a_1,a_2,\dots]$ and $m\ge 1$,
we denote by
\begin{equation} \label{eq:alfam}
\alpha_m := [a_{m+1}, a_{m+2}, \dots]=G^m(\alpha)
\end{equation}
the $m$-th iterate of $\alpha$ under the Gauss map $G:[0,1]\to
[0,1]$ defined by $G(x)=\{1/x\}$ for $x\ne 0$ and $G(0)=0$. We
denote by $p_n^{(m)}$, $q_n^{(m)}$ and $f_n^{(m)}$ the quantities
corresponding to (\ref{eq:rel-approx}) and (\ref{effenne}) for
$\alpha_m$. It holds
\begin{equation} \label{decr}
f_n^{(m)} =\prod_{k=0}^{n}\alpha_{k+m}
\end{equation}
and therefore
\begin{equation} \label{eq:f2f1}
\alpha_{r+m} = \frac{f_r^{(m)}}{f_{r-1}^{(m)}}, \qquad r\geq 0.
\end{equation}
Let moreover $T_\alpha^{(m)} : X\to X$ denote the rotation with
angle $\alpha_m$ (so that $T_\alpha^{(0)} =T_\alpha$) and let
$t_k^{(m)}$, $r_k^{(m)}$, $r_{k_j}^{(m)}$ and $g_j^{(m)}$ be the
corresponding quantities.

\begin{prop} \label{prop:rotren}
The following relations hold for all $m\ge 0$:
\begin{itemize}
\item[(i)] for all $k\ge 0$
$$t_k^{(m)} = \left\{
\begin{array}{ll}
a_{m+1}+1, & {\rm if}\ \{ r_k^{(m)} \alpha_m \}
< f_1^{(m)} = 1- a_{m+1} \alpha_m\\[0.2cm]
a_{m+1},  & {\rm otherwise}
\end{array} \right.$$

\item[(ii)] for all $j\ge 1$
$$g_{j}^{(m)}=r_{k_j}^{(m)} - r_{k_{j-1}}^{(m)} = \left\{
\begin{array}{ll}
q_2^{(m)}, & {\rm if}\ \{ (j-1) \alpha_{m+2} \} < (1-
\alpha_{m+2})\\[0.2cm]
 q_2^{(m)}+q_1^{(m)}, & {\rm otherwise}
\end{array} \right.$$

\item[(iii)] let $(j_h^{(m)})$ be the subsequence such that
$g_{j_h}^{(m)} = q_2^{(m)}+q_1^{(m)}$ for all $h\ge 0$. Then
$j_0^{(m)} = a_{m+3} +1$ and $j_h^{(m)}- j_{h-1}^{(m)}=
t_{h}^{(m+2)}$ for all $h\ge 1$.
\end{itemize}
\end{prop}

\noindent {\sl Proof.} For notational simplicity' sake we show the
results for $m=0$. The general situation is obviously the same.

\noindent Point $(i)$ has been proved above.

\noindent To prove point $(ii)$ we apply the {\it Three Gap
Theorem} (see for example \cite{pytheasfogg}) to the interval
$(0,f_1)$. According to this theorem the possible values of the
gaps $g_j$ between two successive visits of the interval $(0,f_1)$
by the orbit $(j\alpha)$ are given by
$$g_j=r_{k_j} - r_{k_{j-1}} = \left\{
\begin{array}{ll}
q_2+q_1, & \mbox{with frequency  }\ \frac{f_2}{f_1},
\\[0.3cm] q_2, & \mbox{with frequency  }\ 1-\frac{f_2}{f_1} .
\end{array} \right.$$
Let now $\{r_{k_{j-1}}\alpha\}$ be in $(0,f_1)$. We can repeat the same
argument as for $r_{a_2}$ to prove that $g_j=q_2$ if and only if
$\{r_{k_{j-1}}\alpha\}< f_1-f_2$. Indeed we have $\set{r_{k_{j-1}+1} \alpha}
= \set{r_{k_{j-1}} \alpha} +f_0 - f_1$, and more generally we can write
$$\set{r_{k_{j-1}+h} \alpha} = \set{r_{k_{j-1}} \alpha} +f_0 - h
f_1$$
for all $h=1,2,\dots$ such that the r.h.s. remains
non-negative. This certainly happens until $h$ reaches the value
$a_2$, as one readily checks, but for $h=a_2+2$ we have
$\set{r_{k_{j-1}} \alpha} +f_0 - (a_2+2) f_1 = \set{r_{k_{j-1}}
\alpha} + f_2 - 2f_1 < f_2 - f_1<0$ since $\set{r_{k_{j-1}}
\alpha} < f_1$. This shows that $a_2\le k_j - k_{j-1} \le a_2+1$
and it is a constructive proof of what are the possible values of
the gaps $g_j$. Now $k_j - k_{j-1}$ is equal either to $a_2$ or to
$a_2+1$ (and the gap $g_j$ is equal to $q_2$ or to $q_2+q_1$,
respectively) if and only if $\set{r_{k_{j-1}} \alpha} +f_0 - a_2
f_1$ is smaller or greater than $f_1$, respectively. Hence
$g_j=q_2$ if and only if $\{r_{k_{j-1}}\alpha\}< f_1-f_2$.

\noindent Let us now denote by $\tilde T$ the map that acts on
$\tilde X := [0,f_1]/(0=f_1)$ as $\tilde T : \set{r_{k_j} \alpha}
\mapsto \set{r_{k_{j+1}} \alpha}$, that is the first return map on
the interval $[0,f_1]$ for the rotation $T_\alpha$. $\tilde T$ is
isomorphic to the rotation $T_\alpha^{(2)}$ of $X$ through the
angle $\alpha_2= \frac{f_2}{f_1} =[a_3,a_4,\dots]$. Starting from
$r_{k_0}=0$, we need to follow the orbit of $0$ using the rotation
$T_\alpha^{(2)}$ and determine the gaps $(g_j)$. We have showed
that $g_j =q_2$ if and only if $(T_\alpha^{(2)})^{j-1}(0) <
1-\frac{f_2}{f_1}$ or, which is the same by (\ref{eq:f2f1}), $\{
(j-1) \alpha_2 \} < (1-\alpha_2)$. This proves $(ii)$.

\noindent Point (iii) follows by repeating a similar argument for
$T_\alpha^{(2)}$ on the interval $(1-\alpha_2, 1)$. Again the
Three Gap Theorem yields the gaps between two successive visits of
the interval $(1-\alpha_2, 1)$ by the orbit $((j-1)\alpha_2)$. Let
$(j_h)$ be the subsequence such that $\{ (j_h -1) \alpha_2 \} >
(1-\alpha_2)$, then for each $h\ge 1$, ${j_h} - {j_{h-1}}$ is
equal either  to $a_3+1$ or to $a_3$, with ${j_0}= a_3+1$. Indeed,
$(j_0-1) \alpha_2 = a_3 \alpha_2 <1$ but $(a_3+1) \alpha_2 >1$,
therefore $a_3 \alpha_2 > (1-\alpha_2)$. This implies that
$g_1=g_2=\dots=g_{a_3} = q_2$ and $g_{a_3+1}= q_2+q_1$. Let us
remark that $\lfloor (j_{h-1}-1) \alpha_2 \rfloor =h-1$, then
$r_h^{(2)}=j_{h-1}$, where we recall that $r_h^{(2)}$ is defined
as the smallest integer such that $\lfloor r_h^{(2)} \alpha_2
\rfloor =h$. Now $\{ r_h^{(2)} \alpha_2 \} + (a_3+1) \alpha_2 >1$
hence $(j_h -1)-(j_{h-1}-1) \le a_3+1$, and $\{ r_h^{(2)} \alpha_2
\} + (a_3-2) \alpha_2 < (a_3-1)\alpha_2 < 1-\alpha_2$ hence $(j_h
- 1) - (j_{h-1} - 1) \ge a_3$. Moreover $\{ r_h^{(2)} \alpha_2 \}
< f_1^{(2)}=1-a_3 \alpha_2$ if and only if $\{ r_h^{(2)} \alpha_2
\} + (a_3-1) \alpha_2 < 1-\alpha_2$, hence if and only if $(j_h
-1)-(j_{h-1}-1) = a_3+1$. This shows that for all $h\ge 1$ ${j_h}
- {j_{h-1}}$ is equal to $t_h^{(2)}$ and both are equal either to
$a_3+1$ or to $a_3$. $\Box$

\vskip 0.5cm
\noindent The proof given above brings out the {\sl
renormalization argument} mentioned in the Introduction and which
will be fully developed in the next section. According to the
above discussion, since $j_0=a_3+1$,  the values of the sequence
$r_{k_j}$, with $1\le j \le a_3+1$, are given by
$$0,\ q_2,\ 2q_2,\dots,\ a_3 q_2,\ a_3 q_2 + q_2 +q_1 = q_3+q_2.$$
Note that $k_{a_3+1}=p_3 + p_2$. To continue the determination of
the numbers $r_{k_j}$ we have to use the knowledge of the
following $j_h$ and by point (iii) this is equivalent to repeat
the argument above for the rotation $T_\alpha^{(2)}$.

\noindent We need the \emph{Ostrowski representation} of an
integer number \cite{pytheasfogg}: given an irrational number
$\alpha \in (0,1)$ with partial quotients $(a_h)$ and denominators
$(q_h)$ of its rational approximants, any positive integer $r$ can
be written in a unique way in the form
\begin{equation} \label{eq:ost}
r = \sum_{h\ge 0}^N c_h \,q_h \quad \hbox{with}\quad 0\leq c_{h}
\leq a_{h+1} \quad \hbox{and}\quad c_{h-1}=0\quad \hbox{if}\quad
c_{h}=a_{h+1}
\end{equation}
for some integer $N$. We call $N$ the \emph{order} of the integer
$r$, denoted as $N=ord(r)$.

\begin{thm} \label{cor:rkj}
Given a positive integer $r$, we have $r=r_{k_j}$ for some $j> 0$
if and only if in the Ostrowski representation of $r$ we have:
$c_0=c_1=0$ and $\min \{h : c_h >0\}\geq 2$ and even. Moreover, a
positive integer $r$ is of the form $r_k$, for some $k$, if and
only if either $r=r_{k_j}$ for some $j> 0$ or $r=r_{k_j} + c_1 q_1
+1$ for some $j> 0$ and $1\le c_1 \le a_2$.
\end{thm}
\noindent {\sl Proof.} We have verified the first part of the
thesis for $r\le q_3+q_2$, finding $r_{k_j}$ with
$j=1,\dots,a_3+1$.
To continue, by Proposition \ref{prop:rotren}(iii) we
have to study the sequence $t_h^{(2)}$, that is the sequence $t_h$
for the angle $\alpha_2=[a_3,a_4,\dots]$. The first $(a_4+1)$
values are
$$a_3+1,\ a_3,\ a_3,\dots,\ a_3, a_3+1$$ as obtained by part (i)
and (ii) of Proposition \ref{prop:rotren} for $m=2$. This
leads to the computation of $r_{k_j}$ up to $(q_4+q_3+q_2)$. What
happens after depends on whether $r_{k_2}^{(2)}$ is
$q_2^{(2)}=(a_4 a_3+1)$ or $q_2^{(2)}+ q_1^{(2)}=(a_4 a_3 +a_3
+1)$. We already solved this problem for $r_{k_j}^{(0)}$ up to
$j=a_3+1$. Hence in the same way we can solve the problem for
$r_{k_j}^{(2)}$ up to $j=a_3^{(2)}+1=a_5+1$. This implies the
thesis up to $(q_5+q_2)$.

\noindent
The subsequent steps follow by repeating the same argument as
before, where for all $i\ge 2$, the denominators $q_{2i}$ and
$q_{2i+1}$ substitute $q_4$ and $q_5$. Whence the form of the
integers $r_{k_j}$ follows by induction on $i\ge 2$.

\noindent
To prove the result for $r_k$, simply notice that if
$r_k$ is not $r_{k_j}$, then it is obtained from one of the
$r_{k_j}$ by adding $q_1$ as many times as needed, since $t_k\ge
q_1$, hence there are at least $q_1$ iterations before $\lfloor
r\alpha \rfloor$ increases. Moreover the iterations can't be more
than $q_1$ because $r_k$ is not $r_{k_j}$. This finishes the
proof. $\Box$

\vskip 0.5cm \noindent We finally point out that the following
relation between an integer $r_{k_j}$ and its index $j$ is in
force: first, for $j>0$ we have
\begin{equation} \label{eq:rkj-kj}
r_{k_j} = \sum_{h\ge 2} c_h\, q_h \qquad \Rightarrow \qquad k_j=
\sum_{h\ge 2} c_h\, p_h.
\end{equation}
Second, replacing $p_3$ and $p_2$ in $k_j$ with $a_3$ and $1$,
respectively, and using the definition of the numbers
$q_h^{(2)}$, one obtains inductively
\begin{equation} \label{eq:kj-j}
j= \sum_{h\ge 2} c_{h} q_{h-2}^{(2)} = \sum_{h\ge 2} c_{h} \left(
q_h - a_1 p_h \right).
\end{equation}

\noindent
In the following, besides $\alpha_m=G^m(\alpha)$ we will
also need the numbers
\begin{equation} \label{eq:alfamprimo}
\bar \alpha_m := [a_{m+1}-1, a_{m+2},\dots] =
\frac{G^m(\alpha)}{1-G^m(\alpha)}\, \cdot
\end{equation}
We remark that if $a_{m+1}=1$ then $\bar \alpha_m= \alpha_{m+2}$.
Let us denote by $\bar T_\alpha^{(m)}:X\to X$ the rotation of
angle $\bar \alpha_m$ and $\bar p_n^{(m)}, \bar q_n^{(m)}, \bar
f_n^{(m)}$ the corresponding quantities (cfr.
(\ref{eq:rel-approx}) and (\ref{effenne})). A simple inductive
argument shows that in the case $m=1$, if $r_{k_j}$ is defined as
above, then
\begin{equation} \label{eq:kjmenoj}
k_j-j = \left\{
\begin{array}{ll}
\sum_{h\ge 2} \ c_h \bar q_{h-1}^{(1)}, & \mbox{if } a_2 \not= 1, \\[5mm]
\sum_{h\ge 3} \ c_h \bar q_{h-3}^{(1)}, & \mbox{if } a_2 = 1.
\end{array} \right.
\end{equation}
Whereas for the sequence $(\bar f_n^{(m)})$ it holds for all $m\ge
1$ and all $n\ge 0$
\begin{equation} \label{eq:effebar}
\bar f_n^{(m)} = \left\{
\begin{array}{ll}
\frac{f_{n+1}^{(m-1)}}{f_0^{(m-1)}-f_1^{(m-1)}}, & \mbox{if } a_2 \not= 1, \\[5mm]
\frac{f_{n+3}^{(m-1)}}{f_2^{(m-1)}}, & \mbox{if } a_2 = 1.
\end{array} \right.
\end{equation}

\vskip 0.5cm
\noindent We end this section by giving the following
version of a standard expansion of a real number $\beta \in (0,1)$
in terms of the numbers $f_n$ defined in (\ref{effenne}) for a
fixed irrational number $\alpha$ with partial quotients $(a_k)$
(see, e.g., \cite{pytheasfogg}, Sect. 6.4)

\begin{prop} \label{prop:espansione}
For all $\beta \in (0,1)$ there exists a unique sequence of
integers $(b_k)$ such that: (i) $\beta= \sum_{k=0}^\infty b_k
f_k$; (ii) $0\le b_k \le a_{k+1}$ for all $k\ge 0$; (iii)
$b_k=a_{k+1}$ implies $b_{k+1}=0$. Moreover the coefficients
$(b_k)$ are definitively null if and only if $\beta \in \Z +
\alpha \Z$.
\end{prop}

\noindent \emph{Proof.} By definition, $(f_k)$ is a monotonically
decreasing sequence of positive real numbers. The sequence $(b_k)$
is constructed by a greedy algorithm: let
$$b_0 := \left\lfloor \frac \beta \alpha \right\rfloor, \qquad \qquad \beta_1
:= \beta - b_0 \alpha,$$ where we recall $\alpha=f_0$. Note that
$\beta < 1$ implies $b_0 \le a_1$ and $\beta_1 < f_0$. Then we can
define by induction for all $k\ge 1$
\begin{equation} \label{eq:bk}
b_k := \left\lfloor \frac{\beta_k}{f_k} \right\rfloor, \qquad
\qquad \beta_{k+1} := \beta_k - b_k f_k = \beta - \sum_{i=0}^k b_i
f_i.
\end{equation}
By definition of $b_k$, it holds $\beta_k < f_{k-1}$, hence
$\beta=\lim_k \sum_{i=0}^k b_i f_i$ and $b_k \le a_{k+1}$ (see
equation (\ref{recu1})). Moreover, $b_k = a_{k+1}$ implies
$\beta_{k+1} < f_{k-1} - a_{k+1} f_k = f_{k+1}$ (see equation
(\ref{recu})), hence $b_{k+1} =0$. This proves part (i), (ii) and
(iii).

\noindent Let now $b_k=0$ for all $k> \bar k$, for some integer
$\bar k$. Then $\beta = \sum_{i=0}^{\bar k} b_i f_i$ and $f_k \in
\Z + \alpha \Z$ for all $k\ge 0$ imply that $\beta \in \Z + \alpha
\Z$. Conversely, let $\beta = t+\alpha s$ for $t,s \in \Z$. Since
$\beta \in (0,1)$, if $t=0$ then $0\le s \le a_1$, hence $b_0=s$
and $\beta_1 =0$. This implies $b_k=0$ for all $k\ge 1$. Let now
$t>0$ so that $s<0$. If we let $m=\max \set{r \in \N \ :\ \lfloor
r\alpha \rfloor =t-1}$, then we can write $\beta = t -m\alpha +
(m-|s|)\alpha$, where $m-|s| \le a_1$. Using the expansion of
equation (\ref{eq:miaforma}), we can write $m=r_{k_j}+R_1 q_1$ for
some $r_{k_j}=\sum_{h=2}^N c_h q_h$ and $0\le R_1 \le a_2+1$.
Notice that $R_0=0$ by the definition of $m$. From this, using
equation (\ref{eq:rkj-kj}), we obtain $t= k_j +R_1$ and therefore
$$\beta = \sum_{h=2}^N \ (-1)^{h+1} c_h f_h + R_1 f_1 + (m-|s|)
f_0$$ From the definition of $b_k$ one immediately sees that
$b_k=0$ for all $k> N$. The same argument works for the case
$t<0$. $\Box$

\section{The growth of $S_n(\alpha)$} \label{sec:diffusione}

We now use the sequence $(t_k)$ to study the behaviour of
$S_n(\alpha)$ and whence the diffusive properties of the
corresponding walk.

\noindent We have showed that $t_k$ is equal either to $a_1$ or to
$a_1+1$. Therefore, according to whether $a_1$ is even or odd,
only the iterations for which $t_k=a_1+1$ or $t_k=a_1$,
respectively, are important for the growth behaviour.

\noindent Let us consider first of all the case $a_1$ even. In
this case, if $t_k=a_1$ then $\sum_{i=r_k}^{r_k+a_1-1}
s(i\alpha)=0$, hence we can neglect these terms, since the
\virg{walker} associated to $S_n(\alpha)$ would simply take $a_1$
steps to start from $S_{r_k-1}(\alpha)$ and come back to the same
point, after having reached the point
$S_{r_k-1}(\alpha)+\frac{a_1}{2}$. Hence we can restrict ourselves
to the study of the sequence $\set{r\alpha}$ with $r_{k_j}\le r
\le r_{k_j}+a_1$, where we recall that the sub-sequence
$(r_{k_j})$ corresponds to $t_{k_j} = a_1+1$. In these cases
$\sum_{i=r_{k_j}}^{r_{k_j}+a_1} s(i\alpha)= \pm 1$, according to
whether the number $\{ r_{k_j} \alpha + \frac{a_1}{2} \alpha \}$
is $< \frac 1 2$ or $> \frac 1 2$, respectively, that is whether
$\set{r_{k_j} \alpha} < \frac 1 2 f_1$ or $> \frac 1 2 f_1$. In
view of the analysis made in the previous section, given the first
return map $\tilde T$ on the interval $(0,f_1)$, and its
isomorphism with the rotation $T_\alpha^{(2)}$ on $X$, we conclude
that
$$\sum_{i=r_{k_j}}^{r_{k_j}+a_1} s(T_\alpha^i(0))= 1 \
\Longleftrightarrow \ (T_\alpha^{(2)})^j (0) < \frac 1 2\, \cdot$$
Using this fact we now study the relation between the sequences
$S_n(\alpha)$ and $S_n(\alpha_2)$. We obtain that for all $r\ge 0$
it holds
\begin{equation} \label{eq:rin-rotpari}
a_1 \mbox{ even} \ \Longrightarrow \ S_{r}(\alpha) =
S_{j(r)}(\alpha_2) + \tilde S(r)
\end{equation}
where $j(r)$ and $\tilde S(r)$ are computed in the following way.
Let us write $r$ in the form
\begin{equation} \label{eq:miaforma}
r= r_{k_j} + R_1 q_1 +R_0
\end{equation}
with $R_1 q_1 + R_0 < r_{k_{j+1}}- r_{k_j}$, $0\le R_1 \le a_2+1$
and $0\le R_0 < q_1$. We remark that this can be different from
the Ostrowski representation of $r$, since it can be $R_1=a_2+1$.
However the order of $r$ is equal to that of $r_{k_j}$ for $j>0$.
We have
\begin{equation} \label{eq:jr}
j(r)=\max \set{\bar j\ge 0 \ :\ r_{k_{\bar j}} < r-R_0} = \max
\set{j+ \mbox{sgn}(R_1)-1,0}
\end{equation}
and
\begin{equation} \label{eq:sr}
\tilde S(r)= \left\{
\begin{array}{ll}
\sum_{i=r-R_0+1}^r s(\{ i\alpha \}) & \mbox{ if } R_0>0, \; R_1>0, \\[0.3cm]
\sum_{i=r-R_0}^r s(\{ i\alpha \}) & \mbox{ if } R_0>0, \;R_1=0,\\[0.3cm]
0& \mbox{ if } R_0=0, \;R_1>0,\\[0.3cm]
s(\{ r\alpha \})& \mbox{ if } R_0=0, \;R_1=0.
\end{array} \right.
\end{equation}
We remark that using equations (\ref{eq:rkj-kj}) and
(\ref{eq:kj-j}) it is possible to obtain $j$ from the knowledge of
$r_{k_j}$. Moreover $0\le \tilde S(r) \le 1 + \frac{a_1}{2}$ for
all $r\ge 0$, hence the growth behaviour of $S_n(\alpha)$ only
depends on that of $S_n(\alpha_2)$.

\vskip 0.5cm \noindent The case $a_1$ odd is in some sense
complementary to the previous one. Indeed, in this case, we
obviously have $\sum_{i=r_{k_j}}^{r_{k_j}+a_1} s(T_\alpha^i(0))=
0$, whereas for $k$ such that $t_k = a_1$ we have
$\sum_{i=r_{k}}^{r_{k}+a_1-1} s(T_\alpha^i(0))= \pm 1$  according
to whether $\{r_k \alpha + \frac{a_1-1}{2} \alpha \} < \frac 1 2$
or $>\frac 1 2$. We would like to construct an induced map on some
interval of $X$, to connect the values of
$\sum_{i=r_{k}}^{r_{k}+a_1-1} s(T_\alpha^i(0))$ to a suitable
orbit of such induced map. To this aim we notice that the point
$\{ r_k \alpha + \frac{a_1-1}{2} \alpha \}$ belongs to the
interval $J:=\left( \frac 1 2 - \frac 1 2 ( f_0 - f_1), \frac 1 2
+ \frac 1 2 ( f_0 - f_1)\right)$ for all $k\ge 0$ such that $t_k
=a_1$. This follows immediately from the following remarks:
\begin{enumerate}
\item
$\set{r_k \alpha}> f_1$ and $f_1 + \frac{a_1-1}{2} \alpha = \frac
1 2 - \frac 1 2 ( f_0 - f_1)$;

\item $\{ r_k \alpha +
\frac{a_1-1}{2} \alpha \} < f_1 + \frac{a_1+1}{2} \alpha - f_1 $,
since $f_1 + \frac{a_1+1}{2} \alpha > \{ r_{k_j} \alpha +
\frac{a_1-1}{2} \alpha \}$ for all $k_j$, and $\frac{a_1+1}{2}
\alpha = \frac 1 2 + \frac 1 2 ( f_0 - f_1)$.

\end{enumerate}
Moreover the two
estimates in 1. and 2. are sharp.

\noindent From the definition of the interval $J$, it also follows
that $\{ r_{k_j} \alpha + r \alpha \} \not\in J$ for all
$r=1,\dots,a_1$, since $\set{r_{k_j} \alpha} \in (0,f_1)$ implies
that $\{ r_{k_j} \alpha + \frac{a_1-1}{2} \alpha \} <\frac 1 2 -
\frac 1 2 ( f_0 - f_1)$ and $\{ r_{k_j} \alpha + \frac{a_1+1}{2}
\alpha \} > \frac 1 2 + \frac 1 2 ( f_0 - f_1)$. Hence we can
consider the first return map $\bar T$ of $T_\alpha$ to the
interval $J$, and obtain that $\bar T$ is isomorphic to the
inverse of the rotation $\bar T_{\alpha}^{(1)}$ on $X$, that is
the rotation of angle $$- \bar \alpha_1=- \frac{f_1}{f_0-f_1}.$$

\noindent
Let now $(k_i)$ be the sub-sequence such that
$t_{k_i}=a_1$ for all $i\ge 1$, then
$$\sum_{r=r_{k_i}}^{r_{k_i}+a_1-1} s(T_\alpha^r(0))= 1
\Longleftrightarrow \ (\bar T_\alpha^{(1)})^i (0) < \frac 1 2$$ We
now want to give an analogous equation of (\ref{eq:rin-rotpari}).
In this case we have to neglect $(\bar T_\alpha^{(1)})^0(0)=0$,
since we start with $i=1$, hence we obtain that for all $r\ge 0$
\begin{equation} \label{eq:rin-rotdisp}
a_1 \mbox{ odd} \ \Longrightarrow \ S_{r}(\alpha) =
S_{i(r)}(-\bar \alpha_1)-1 + \tilde S(r)
\end{equation}
where $\tilde S(r)$ is the same as in equation
(\ref{eq:rin-rotpari}) and $i(r)$ is computed in the following
way. Let us write again $r$ as in equation (\ref{eq:miaforma}),
then
\begin{equation} \label{eq:ir}
i(r)= k_j-j+\max\{(R_1-1),0\}
\end{equation}
where we recall equations (\ref{eq:rkj-kj}), (\ref{eq:kj-j}) and
(\ref{eq:kjmenoj}).

\noindent
Again $\tilde S(r)$ is uniformly bounded so that the diffusive properties of
$S_n(\alpha)$ depend only on those of
$S_n(\bar \alpha)$. Moreover we note that for all $n\ge 0$ we have
\begin{equation} \label{eq:utile-dispari}
S_{n}(-\bar \alpha_1) -1  = - \left( S_n(\bar \alpha_1)-1 \right)
\end{equation}

\vskip 0.5cm \noindent In conclusion, we have showed that, as far
as the diffusive properties are concerned, the walk
$(S_n(\alpha))$ is equivalent to a \virg{renormalized} walk
$(S_{R(n)}(\beta))$, where the values of $R(n)$ and $\beta$ depend
on the parity of $a_1$, the first partial quotient of the number
$\alpha$.

\noindent Equations (\ref{eq:rin-rotpari}) and
(\ref{eq:rin-rotdisp}) lead by iteration to an explicit expression
for $S_n(\alpha)$ only in terms of the $(a_k)$. A tentative result
in this direction was given in \cite{sos}. About growth estimates,
let us see how this argument leads to precise estimates on the
behaviour of maxima and minima of $S_n(\alpha)$.

\begin{prop} \label{prop:massimi}
Given $\alpha=[a_1,a_2,\dots] \in (0,1)$, let $r=r_{k_j}+R_1 q_1 +
R_0$ for some $j\ge 0$ as in equation (\ref{eq:miaforma}).
\noindent If $a_1$ is even then
$$0 \le \max\limits_{0\le n\le r} S_n(\alpha) - \left( \max\limits_{0\le m\le j(r)}
S_m(\alpha_2) +\frac{a_1}{2} \right) \le 1$$
$$\min\limits_{0\le n\le r} S_n(\alpha) = \min\limits_{0\le m\le j(r)}
S_m(\alpha_2)$$ where $\alpha_2=[a_3,a_4,\dots]$ and $j(r)$ is given in (\ref{eq:jr}). If instead
$a_1$ is odd then
$$0\le \max\limits_{0\le n\le r} S_n(\alpha) - \left(1-\min\limits_{0\le m\le i(r)}
S_m(\bar \alpha_1) + \frac{a_1-1}{2}\right) \le 1$$
$$\min\limits_{0\le n\le r} S_n(\alpha) = 1-\max\limits_{0\le m\le i(r)}
S_m(\bar \alpha_1)$$ where $\bar \alpha_1 = [a_2-1,a_3,\dots]$
and $i(r)$ is given in (\ref{eq:rin-rotdisp}).

\noindent Moreover in the case $a_1$ even, the difference between
maxima is equal to 1 only if $R_1=0$ and $R_0\ge \frac{a_1}{2}$.
\end{prop}

\noindent \emph{Proof.} Let us consider first the case $a_1$ even.
The result is a direct consequence of equation
(\ref{eq:rin-rotpari}) and the relation $0\le \tilde S(r) \le 1 +
\frac{a_1}{2}$.

\noindent
For the case $a_1$ odd, the proof follows from equations
(\ref{eq:rin-rotdisp}) and (\ref{eq:utile-dispari}). $\Box$

\vskip 0.5cm
\noindent
We point out that, since $j(r)$ and $i(r)$ are explicitly
computable from $r$, one can iterate the renormalization argument
in such a way that at the end of the process the maxima and minima of the
walk $S_n(\alpha)$ will be explicitly computable linear
combinations of the partial quotients of $\alpha$. To this end we observe that in order to apply
the argument to $S_{m}(\alpha_2)$ it
is enough to notice that in equation (\ref{eq:kj-j}) the number
$j$ is obtained as a linear combination of $a_3$ and $1$, which are
nothing but $q_1^{(2)}$ and $q_0^{(2)}$ respectively (we are using
the notations of Proposition \ref{prop:rotren}). Therefore $j$ can
be expressed with respect to $\alpha_2$ in the form $r_{k_l}^{(2)}
+ R_1^{(2)} q_1^{(2)} + R_0^{(2)}$, and the iteration can proceed.
We thus obtain the following,

\begin{thm} \label{cor:maxpari}
Let the partial quotients $(a_{2i+1})$ be even for all $i\ge 0$.
If $r=r_{k_j}+R_1 q_1 + R_0$ for some $j\ge 0$ and $ord(r)=N$,
then
$$\frac{1}{2} \sum_{i=0}^{\frac{N-2}{2}} a_{2i+1} \le
\max\limits_{0\le n\le r} S_n(\alpha) \le \frac N 2 +
\frac{1}{2} \sum_{i=0}^{\frac{N-2}{2}} a_{2i+1}$$
$$\min\limits_{0\le n\le r} S_n(\alpha) =1$$
\end{thm}
\noindent This theorem implies that the diffusion properties of
$S_n(\alpha)$ depend only weakly on the partial quotients
$(a_{2i})$. In particular, for all $\alpha$ with fixed partial
quotients $(a_{2i+1})$, even for all $i\ge 0$, the sequence
$S_n(\alpha)$ grows with the same rate, and what changes is the
number of fluctuations.

\noindent The situation is more cumbersome for numbers $\alpha$
with odd partial quotients in an odd position. This would imply to
change the kind of \virg{renormalization}, and also partial
quotients with even position become important. However, we can
make some computations for particular cases.

\vskip 0.2cm \noindent {\sc Example.} Let $\alpha=[a,a,a,\dots]$
with $a$ odd. Then the first renormalization leads to $\bar
\alpha_1 =[a-1,a,a,\dots]$. This fact implies that two different
situations occur for $a=1$ and $a>1$. Hence the sequence
$S_n(\alpha)$ with $\alpha$ the golden ratio
$\frac{\sqrt{5}-1}{2}$ has peculiar properties.

\noindent
Let us first consider the case $a>1$. From Proposition
\ref{prop:massimi} it follows that
$$\max\limits_{0\le n\le r} S_n(\alpha) \le 2 + \frac{a-1}{2}
- \min\limits_{0\le m\le j(i(r))} S_m(\alpha) \le $$
$$\le 2+ (a-1) + \max\limits_{0\le n\le j(i(j(i(r))))}
S_n(\alpha)$$ where $ord(j(i(j(i(r)))))=ord(r)-6$. Therefore, if
$ord(r)=6k$, repeating the same argument from below we have
$$\frac{(a-1)}{6}\ ord(r)\le \max\limits_{0\le n\le r} S_n(\alpha)
 \le \frac{(a+1)}{6}\ ord(r)$$
For example, if $a=3$ then $\alpha=\frac{\sqrt{13}-3}{2}$, and
$$\frac{1}{3\ \log(\frac{\sqrt{13}+3}{2})}\le \limsup\limits_{r \to \infty}\ \frac{
\max\limits_{0\le n\le r} S_n(\frac{\sqrt{13}-3}{2})} {\log r}
\le \frac{2}{3\ \log(\frac{\sqrt{13}+3}{2})}$$

Let us consider now $a=1$. From Proposition \ref{prop:massimi} it
follows that
$$\max\limits_{0\le n\le r} S_n\left( \frac{\sqrt{5}-1}{2}\right) \le 1+
\max\limits_{0\le m\le i(i(r))}
S_m\left(\frac{\sqrt{5}-1}{2}\right)$$ and
$ord(i(i(r)))=ord(r)-6$. Hence
$$\limsup\limits_{r \to \infty}\ \frac{
\max\limits_{0\le n\le r} S_n\left(\frac{\sqrt{5}-1}{2}\right)}
{\log r} \le \frac{1}{6\ \log(\frac{\sqrt{5}+1}{2})}.$$

\vskip 0.5cm

%\cgraf(goldenratio.eps scaled 1000) Tit= Testo=$S_n(\alpha)$ for
%$\alpha= [1,1,1,\dots]$ with $0\leq n\leq 2000$ X= Y={}

\vskip 0.5cm

\section{Generalisations to other orbits and applications}
\label{sec:disc}

We now use the formalism developed in the previous sections to
analyse the behaviour of the following quantity
\begin{equation} \label{def:altri-punti}
d_n(\alpha,\beta) := \sum_{r=0}^{n-1} \ \chi_{[0,\beta)}
(\set{r\alpha}) - \beta n
\end{equation}
for points $\beta \in (0,1)$. We call $d_n(\alpha,\beta)$ the
\emph{relative discrepancy} of $\alpha$ with respect to $\beta$.
The term is justified by the usual definition of
\emph{discrepancy} of the sequence $(n\alpha)$ as
\begin{equation} \label{def:disc}
D_n^*(\alpha):= \sup\limits_{\beta \in (0,1)} \left| \frac 1 n \
d_n(\alpha,\beta) \right|
\end{equation}
We first give an iterative argument to compute the relative
discrepancies $d_n(\alpha,\beta)$. This is useful to give an
explicit expression for the diffusion of the orbit of a general
point $\beta\in (0,1)$ for a rotation $T_\alpha$, that is
$$S_n(\alpha,\beta):= \sum_{r=0}^n s(\set{r\alpha + \beta})$$
(see (\ref{def:somme})).

\noindent
Given $\alpha=[a_1,a_2,\dots]$, we recall the notations
$$\alpha_m=G^m(\alpha)=[a_{m+1},a_{m+2},\dots]$$
$$\bar \alpha_m=\frac{G^m(\alpha)}{1-G^m(\alpha)}=[a_{m+1}-1,a_{m+2},\dots]$$
where $G$ is the Gauss map, as well as the sequences $p_h$, $q_h$, $f_h$
and $r_{k_j}$ associated to $\alpha$. Let us fix $\beta \in
(0,1)$ and recall its expansion $\beta=\sum_k b_k f_k$ as well as the numbers $\beta_m$ given
in (\ref{eq:bk}). Let us define for all $m\ge 1$
\begin{equation} \label{eq:betabar}
\beta^m := \frac{\beta_m}{f_m} \qquad \qquad \bar \beta^m:=
\frac{\beta_m-f_m}{f_{m-1}-f_m}
\end{equation}

\begin{prop} \label{prop:discr1}
For a given $n\in \N$, let us write $n-1=r_{k_j}+R_1 q_1 + R_0$ as
in equation (\ref{eq:miaforma}) with $ord(n)=N$. Let $j(n-1)$ and
$i(n-1)$ be defined as in equations (\ref{eq:jr}) and
(\ref{eq:ir}) respectively, and let $\tilde S(n-1)$ be defined as
in equation (\ref{eq:sr}). If we define
$$S(n,\alpha,\beta)=\tilde S(n-1)-\beta (R_0+1 - \mbox{\rm sgn}(R_1))$$
and write $r_{k_j}=\sum_{h=2}^N c_h q_h$, then
$$d_n(\alpha,\beta) = S(n,\alpha,\beta) + \left\{
\begin{array}{l}
-\frac{(b_0-a_1 \beta)}{f_1} \ C(\alpha,n) +
d_{j(n-1)+1}(\alpha_2,\beta^1)  \\[0.5cm]
\frac{(b_0-a_1 \beta+1-\beta)}{f_0-f_1} \ C(\alpha,n) + \bar
\beta^1 + d_{i(n-1)+1}^c(\bar \alpha_1,\bar \beta^1)
\end{array} \right.$$ where the first formulation is valid if
$b_1=0$ and the second otherwise, the constant $C(\alpha,n)$ does
not depend on $\beta$ and is given by
$$C(\alpha,n):=\alpha\ \mbox{\rm
sgn}(R_1)- R_1 f_1 + \sum_{h=2}^N \ (-1)^{h} c_h f_h $$ and
$d^c(\cdot)$ means that in the definition of $d(\cdot)$ we use the
indicator function of the interval $(1-\beta^{1},1]$.
\end{prop}

\noindent \emph{Proof.} The main idea is to use the partition of
the sequence $\set{r\alpha}$ using the sequence $\set{r_k}$.
Indeed, as in the treatment of the diffusion, we have that for all
$k\ge 0$
$$\sum_{r=r_k}^{r_{k+1}-1} \ \chi_{[0,\beta)}(\set{r\alpha}) =
\left\{
\begin{array}{ll}
b_0+1 & \mbox{if }\ \{r_k \alpha\} < \beta_1 \\[0.2cm]
b_0 & \mbox{otherwise}
\end{array} \right.$$
By lemma \ref{lemmino}, if $b_1=0$, that is $\beta_1<f_1$, the first case is
possible only if $k=k_j$ for some $j\ge 0$. In this case we have
$$\sum_{r=0}^{n-1} \ \chi_{[0,\beta)}(\set{r\alpha})=\tilde S(n-1)+
b_0(k_j+R_1)+\sum_{j=0}^{j(n-1)} \
\chi_{[0,\beta^1)}(\set{j\alpha_2})$$ where the last term accounts
for the relation
$$\sum_{r=r_{k_j}}^{r_{k_j}+q_1} \
\chi_{[0,\beta)}(\set{r\alpha})= b_0 +
\chi_{[0,\beta_1)}(\set{r_{k_j}\alpha})$$ and uses the isomorphism
between the first return function to the set $(0,f_1)$ and the
rotation of angle $\alpha_2$ on the unit circle. Moreover let us
write
$$\beta n = \beta (R_0+1 - \mbox{\rm sgn}(R_1)) + \beta (j+\mbox{\rm
sgn}(R_1)) + \beta (r_{k_j}+R_1q_1-j)$$ where, we recall,
$j(n-1)=j+\mbox{\rm sgn}(R_1)-1$ (cfr. (\ref{eq:jr})) and
$(r_{k_j}+R_1q_1-j)=a_1(k_j+R_1)$ (cfr. (\ref{eq:rkj-kj})
and (\ref{eq:kj-j})). The claim now follows by evaluating
$$k_j+R_1
- \frac{\alpha}{f_1} (j+\mbox{\rm
sgn}(R_1))=-\frac{C(\alpha,n)}{f_1}.$$ A similar argument can be
applied if $b_1>0$, using the fact that the first return map to
the interval $(f_1, \alpha)$ is isomorphic to the inverse of the
rotation on the unit interval with angle $\bar \alpha_1$.
This leads to  the term $d^c(\cdot)$. $\Box$

\vskip 0.5cm
\noindent
The previous result
shows that  to evaluate $d_n(\alpha,\beta)$ we have to repeat the same argument for
$d_{j(n-1)+1}$ and $d_{i(n-1)+1}$, respectively. To this end we point out that the same
argument as before yields for general
$\alpha$ and $\beta$
\begin{equation} \label{eq:disc-inv}
d_n^c(\alpha,\beta)= S^c(n,\alpha,\beta) + \left\{
\begin{array}{l}
-\frac{(b_0-a_1 \beta)}{f_1} \ C(\alpha,n) +
d_{j(n-1)+1}^c(\alpha_2,\beta^1)  \\[0.5cm]
\frac{(b_0-a_1 \beta+1-\beta)}{f_0-f_1} \ C(\alpha,n) + \bar
\beta^1 + d_{i(n-1)+1}(\bar \alpha_1,\bar \beta^1)
\end{array} \right.
\end{equation}
where $S^c(n,\alpha,\beta)$ is obtained by using the indicator
function of the interval $(1-\beta,1]$ in $\tilde S(n-1)$.

\noindent
Going on, we see that we can repeat the same argument until the
\virg{renormalized} rotation that we obtain does not have enough
iterates to be renormalized again. When this happens is up to the
order of $n$ and to $\beta$, which implies the
\virg{renormalization path} that have to be followed. Let us see
how the algorithm to choose the new angle of rotation and
the new interval works. We have already seen that the first step is
$$(\alpha,\beta) \longrightarrow \left\{
\begin{array}{ll}
(\alpha_2,\beta^1) & \mbox{ if  } b_1=0 \\[0.2cm]
(\bar \alpha_1, \bar \beta^1)^c & \mbox{ if  } b_1>0
\end{array} \right.$$
and that two subsequent $(\cdot)^c$ cancel out (since they are
generated by two inversion in the rotations). By straightforward
computation, one can readily verify that the general scheme is as
follows: the starting point is always of the form
$(\alpha_m,\beta^{m-1})$ or $(\bar \alpha_m,\bar \beta^m)$, and
 for all $m\ge 2$ we have
\begin{equation} \label{eq:schema}
(\alpha_m,\beta^{m-1}) \longrightarrow \left\{
\begin{array}{ll}
(\alpha_{m+2},\beta^{m+1}), & \mbox{ if  } b_{m+1}=0, \\[0.2cm]
(\bar \alpha_{m+1}, \bar \beta^{m+1})^c, & \mbox{ if  } b_{m+1}>0,
\end{array} \right.
\end{equation}
the same holding true for $(\bar \alpha_m,\bar \beta^m)$.

\noindent To conclude, we remark that the constants
$C(\alpha_m,n)$ and $C(\bar \alpha_m, n)$ are the same as in
Proposition \ref{prop:discr1}, with the values $f_n^{(m)}$ and
$\bar f_n^{(m)}$ computed using equation (\ref{decr}) and
(\ref{eq:effebar}), respectively. Furthermore, the coefficients
multiplying $C(\alpha_m,n)$ and $C(\bar \alpha_m,n)$ at each step
satisfy the following

\begin{lem} \label{lem:disc-limiti}
For all $\alpha$ and $\beta$ it holds:
\begin{itemize}
    \item[(i)] if $b_{m+1}=0$ then $-a_1^{(m)} \le \frac{(b_0^{m}-a_1^{(m)}
    \beta^{m-1})}{f_1^{(m)}} \le a_1^{(m)}$;
    \item[(ii)] if $b_{m+1}>0$ then $-a_1^{(m)} \le \frac{(b_0^{m}-a_1^{(m)} \beta^{m-1}+ 1-
    \beta^{m-1})}{{f_0^{(m)}-f_1^{(m)}}} \le a_1^{(m)}$;
    \item[(iii)] for all $m\ge 0$ and all $n$, it holds
    $0<C(\alpha_m,n)< \alpha_m$.
\end{itemize}
The same relations hold if we consider the corresponding quantities
for $\bar \alpha_m$ and $\bar \beta^{m}$.
\end{lem}
\noindent Using the above algorithm we are able to obtain the
actual values of $d_n(\alpha,\beta)$. To obtain growth estimates
we write
$$d_n(\alpha,\beta)=\C(n,\alpha,\beta)+\res(n,\alpha,\beta)+
\B(n,\alpha,\beta)$$ and give estimates for these three terms
separately. These terms come from the expression of
$d_n(\alpha,\beta)$ given in Proposition \ref{prop:discr1}. The
first term arises by summing the sequence of constants
$C(\alpha_m,n(m))$. The second comes from summation of the terms
$S(n,\alpha,\beta)$. The last term arises by adding the different
$\bar \beta^m$ that we encounter when $b_{m}>0$. In the appendix
we give the proof of the following estimates

\begin{prop} \label{prop:disc-stime}
Let $n$ be an integer such that $ord(n-1)=N$, then for all
$\alpha=[a_1,a_2,\dots]$ and all $\beta$, the following hold:
\begin{itemize}
    \item[(i)] $|\C(n,\alpha,\beta)|< N$;
    \item[(ii)] $|\B(n,\alpha,\beta)| < N$;
    \item[(iii)] $|\res(n,\alpha,\beta)| \le \sum_{m=1}^{N+1}\
    \left( 1+\frac{a_m}{4} \right)$.
\end{itemize}
\end{prop}

\noindent We now compute the sum $S_n(\alpha,\beta)$ for a given
$\alpha$. Let
\begin{equation}
S_n(\alpha,\beta) = S_n(\alpha) + R_n(\alpha,\beta)
\end{equation}
where the term $R_n(\alpha,\beta)$ accounts for the times that
$s(\set{r\alpha}) \not= s(\set{r\alpha +\beta})$. We first remark
that for all $\beta \in [0,\frac 1 2)$ it holds
$$S_n(\alpha, \beta +\frac 1 2) = - S_n(\alpha,\beta)$$
hence it suffices to study $R_n(\alpha,\beta)$ in the case $\beta
< \frac 1 2$.

\begin{thm} \label{prop:gen-diff}
For all $\beta \in [0,\frac 1 2)$,
$$S_n(\alpha,\beta) = S_n(\alpha) + R_n(\alpha,\beta)$$ where
the term $R_n(\alpha,\beta)$ can be written as
$$R_{n-1}(\alpha,\beta) =
2 \left[ d_n(\alpha, \frac 1 2 -\beta) - d_n(\alpha,1-\beta) -
d_n(\alpha,\frac 1 2) \right]$$
\end{thm}

\noindent {\bf Proof.} If we denote
$$P_n := \set{0\le r \le n-1  : \frac 1 2 < \set{r\alpha} < 1
,\ \ 1 < \beta + \set{r\alpha} < \frac 3 2 }$$
$$M_n := \set{0\le r \le n-1  : \set{r\alpha} < \frac 1 2
,\ \ \frac 1 2 < \beta + \set{r\alpha} < 1 }$$ then
$$R_{n-1}(\alpha,\beta) = 2(card(P_n) - card(M_n))$$
Introducing the notations
$$A_a^b:=\set{0\le r \le n-1  : a < \set{r\alpha} < b}$$
$$B_a^b:=\set{0\le r \le n-1  : a < \beta + \set{r\alpha} < b}$$
we have
$$P_n = B_1^{3/2} \cap A_{1/2}^1$$
$$M_n = B_{1/2}^1 \cap A_0^{1/2}= (B_{1/2}^\infty \cap A_0^{1/2})
\setminus (B_1^\infty \cap A_0^{1/2})$$ Now, since $0\le \beta <
\frac 1 2$, it is easy to obtain
$$B_1^{3/2} \cap A^1_{1/2} = B_1^{3/2}$$
$$B_1^\infty \cap A_0^{1/2}= \emptyset$$
hence
$$P_n =B_1^{3/2}=A_{(1-\beta)}^1 = A_0^1 \setminus A_0^{(1-\beta)}$$
$$M_n = B_{1/2}^\infty \cap A_0^{1/2} = A_0^{1/2} \setminus A_0^{(1/2 -\beta)}$$
The thesis now follows by writing
$$card(A_a^b)= \sum_{r=0}^{n-1} \chi_{_{(a,b)}}(\set{r\alpha})$$
and using the definition of $d_n(\alpha,\beta)$ in
(\ref{def:altri-punti}). $\Box$

\vskip 0.5cm \noindent Putting together Proposition
\ref{prop:disc-stime} and Theorem \ref{prop:gen-diff}, one gets
immediately

\begin{cor} \label{cor:l-inf}
For all irrational $\alpha=[a_1,a_2,\dots]$ it holds
$$\sup_\beta |S_n(\alpha,\beta)| \le |S_n(\alpha)| + 6 \sum_{m=1}^N\
\left(3+\frac{a_m}{4} \right)$$ where $ord(n)=N$.
\end{cor}

\subsection{Applications} \label{subsec:disc}

We now give some applications of the estimates of Proposition
\ref{prop:disc-stime}.

We first consider the speed of convergence in the Birkhoff Ergodic
Theorem. As stated in the Introduction, the Lebesgue measure is
invariant and ergodic for the irrational translations on the unit
circle. Hence for the rotation $T_\alpha(\beta):=
\set{\alpha+\beta}$, the Birkhoff Ergodic Theorem implies that
\begin{equation} \label{def:birkhoff}
\lim\limits_{n\to \infty} \frac 1 n \ \sum_{k=0}^{n-1}\ \chi_{I}
(\set{k\alpha +\beta}) = |I|:= \int_0^1\ \chi_{I} (x) \ dx
\end{equation}
for any interval $I\subset [0,1]$. We prove that

\begin{thm} \label{teo:erg}
For any interval $I\subset [0,1]$ it holds
\begin{equation} \label{eq:erg}
\sum_{k=0}^{n-1}\ \chi_{I} (\set{k\alpha +\beta}) - n|I| =
d_n(\alpha,\set{\delta-\beta})-d_n(\alpha,\set{\gamma-\beta})
\end{equation}
from which it follows that
\begin{equation} \label{eq:erg2}
\left|\sum_{k=0}^{n-1}\ \chi_{I} (\set{k\alpha +\beta}) - n|I|
\right| \le 2\ \sum_{m=1}^N\ \left(3+\frac{a_m}{4} \right)
\end{equation}
where $\alpha=[a_1,a_2,\dots]$ and $ord(n-1)=N$.
\end{thm}

\noindent \emph{Proof.} Let $I=[\gamma,\delta]$. Then we can write
$$\sum_{k=0}^{n-1}\ \chi_{I} (\set{k\alpha +\beta}) =
\sum_{k=0}^{n-1}\ \chi_{[\gamma-\beta,\delta-\beta]}
(\set{k\alpha})$$ By using equation (\ref{def:altri-punti}), we
obtain (\ref{eq:erg}). The estimate (\ref{eq:erg2}) follows by
Proposition \ref{prop:disc-stime}. $\Box$

\vskip 0.5cm \noindent One can easily generalize (\ref{eq:erg}) to
bounded variation functions, for which the analogous of estimate
(\ref{eq:erg2}) is the Denjoy-Koksma inequality (see \cite{kn}).
Indeed estimate (\ref{eq:erg2}) can be considered a particular
case of the Denjoy-Koksma inequality.

\vskip 0.5cm

\noindent As stated at the beginning of Section \ref{sec:disc},
the term $d_n(\alpha,\beta)$ is related to the \emph{discrepancy}
$D_n^*(\alpha)$ of the sequence $(n\alpha)$ by equation
(\ref{def:disc}). Some classical results are known for the
discrepancy of a general sequence, and in particular for the
sequence $(n\alpha)$ according to the arithmetical properties of
$\alpha$ (see \cite{kn}). For more recent sharp results we refer
to \cite{sch} and \cite{pinner}. We briefly show below how to get
similar results to those in \cite{sch} and \cite{pinner} directly
by our approach. Proofs can be found in the appendix.

\noindent From Proposition \ref{prop:disc-stime} one gets
immediately
\begin{equation} \label{eq:come-pinner}
nD^*_n(\alpha)=\sup_\beta |d_n(\alpha,\beta)| \le 1+ 3N +\frac 1 4
\ \sum_{m=1}^{N+1} \ a_m
\end{equation}
where $N=ord(n-1)$. This bound is of the same order of results in
\cite{pinner}. We will show that this is the best possible
estimate for the general case. But first we shall obtain some
lower bounds.

\begin{prop} \label{prop:discr2}
For all irrational $\alpha$ there exist a number $\beta$ and two
infinite subsequences $(n_k)$ and $(n_h)$ such that
$$\sup\limits_\beta |\res(n_k,\alpha,\beta)| \ge
\sum_{\stackrel{m=1}{a_{2m-1}\ge 3}}^{\frac{N_k}{2}+1}
\frac{a_{2m-1}-2}{4}$$ where $N_k=ord(n_k)$ and
$$\sup\limits_\beta |\res(n_h,\alpha,\beta)| \ge
\sum_{\stackrel{m=2}{a_{2m}\ge 3}}^{\frac{N_h+1}{2}}
\frac{a_{2m}-2}{4}$$ where $N_h=ord(n_h)$.
\end{prop}

\begin{thm} \label{teo:disc}
Let $\alpha$ have unbounded partial quotients $(a_k)$ and denote
$$\ell_e=\liminf\limits_{k\to \infty} \frac{\sum^k_{m=1\ even}\ a_m}{\sum_{m=1}^k\
a_m}$$
$$\ell_o=\liminf\limits_{k\to \infty} \frac{\sum^k_{m=1\ odd}\ a_m}{\sum_{m=1}^k\
a_m}$$ If $(\ell_e^2 + \ell_o^2) >0$ and
\begin{equation} \label{eq:serve}
\limsup\limits_{k\to \infty}\ \frac{k}{\sum_{m=1}^k\ a_m}=0
\end{equation}
then
\begin{equation} \label{eq:risultato}
\frac{1}{4}\ \max \set{\ell_e, \ell_o}\ \le\ \limsup\limits_{n\to
\infty}\
\frac{nD_n^*(\alpha)}{\stackrel{ord(n-1)+1}{\sum\limits_{m=1}}
a_m} \le \frac 1 4
\end{equation}
\end{thm}

\vskip 0.5cm \noindent {\sc Example.} The conditions of Theorem
\ref{teo:disc} are satisfied by $$\alpha=e-2 =[1,2,1,1,4,1,1,6,1,
1,8,1,1,10,1,1,12,\dots]$$ indeed $\ell_e=\ell_o=\frac 1 2$ and
$\sum_{m=1}^k a_m \sim \frac 1 9 k^2$.

\noindent In fact estimate (\ref{eq:risultato}) is the best
possible in general, as can be shown by choosing for example
$$\alpha:=[1,2,1,3,1,4,1,5,\dots,1,n,\dots]$$
for which $\ell_e=1$, $\ell_o=0$ and $\sum_{m=1}^k a_m \sim \frac
1 8 k^2$. $\lozenge$

\vskip 0.5cm \noindent We finish with few more remarks. First of
all, we remark that by Proposition \ref{prop:espansione} one can
easily prove the well known result that for a given $\alpha \in
(0,1)$, if $\beta \in \Z +\alpha \Z$ then $d_n(\alpha,\beta)$ is
bounded for all $n\ge 0$.

\noindent Finally, we dwell upon relations between the discrepancy
and the sums $S_n(\alpha)$. To start with, let us notice that by
definition $S_{n-1}(\alpha) = 2 d_n(\alpha,\frac 1 2)$, so that
for any given function $F(n)\nearrow \infty$ we have
\begin{equation} \label{eq:disc-somme}
\limsup\limits_{n\to \infty} \frac{nD_n^*(\alpha)}{F(n)} \ge
\limsup\limits_{n\to \infty} \frac{|S_{n-1}(\alpha)|}{2F(n)}
\end{equation}
The analogous relation for the infimum limit is not interesting
since for all $\alpha$ it holds $|S_{q_k}(\alpha)| \le 2$ for all
denominators $q_k$ (see \cite{isola}).

\noindent However, in some cases one could get an equality in
(\ref{eq:disc-somme}). Let us consider the case
$\alpha=\sqrt{2}-1=[2,2,2,\dots]$. We find from Theorem
\ref{cor:maxpari}
$$\limsup\limits_{n\to \infty} \frac{nD_n^*(\sqrt{2}-1)}{\log n}
\ge \limsup\limits_{n\to \infty} \frac{ord(n-1)}{4 \log n} \ge
\frac{1}{4\ \log(\sqrt{2}+1)}$$ This relation is in fact an
equality, as proved in \cite{ds}, and therefore
$$\limsup\limits_{n\to \infty} \frac{S_n(\sqrt{2}-1)}{2 \log n}=
\limsup\limits_{n\to \infty} \frac{nD_n^*(\sqrt{2}-1)}{\log n}$$
The same is shown in \cite{boris2} for
$\alpha=\frac{\sqrt{3}-1}{2}=[2,1,2,1,\dots]$. Moreover the author
exhibits some other couples $(\alpha,\beta)$ for which a similar
relation holds. He then conjectures that a similar relation holds
for all couples $(\alpha,\beta)$ with $\alpha$ a quadratic
irrational and $\beta \in \Q(\alpha)$ but $\beta \not\in \Z +
\alpha \Z$. We show here that it is not the case.

\begin{cor} \label{cor:boris}
Let $\alpha$ be a quadratic irrational in $(0,1)$ with partial
quotients $(a_m)$ verifying $a_{2i-1}=2$ and $a_{2i}=2k$ for all
$i\ge 1$ and a fixed $k>2$. Then $$\nu^*(\alpha) >
\limsup\limits_{n\to \infty} \frac{|d_n(\alpha,\frac 1 2)|}{\log
n}$$
\end{cor}

\noindent \emph{Proof.} Applying the result of \cite{sch}
mentioned above it follows that
$$\limsup\limits_{n\to \infty}
\frac{nD_n^*(\alpha)}{\log n}= \frac 1 4\ \max
\set{\limsup\limits_{N\to \infty}\ \frac{kN}{\log q_N},\
\limsup\limits_{N\to \infty}\ \frac{N}{\log q_N}}$$ However from
Theorem \ref{cor:maxpari} and equation (\ref{eq:disc-somme}), it
follows that if $N=ord(n)$ then
$$\limsup\limits_{n\to \infty}\ \frac{|d_n(\alpha,\frac 1 2)|}{\log
n}\le \limsup\limits_{n\to \infty}\ \frac{N}{2\ \log n} \le
\limsup\limits_{n\to \infty}\ \frac{N}{2\ \log q_N}$$ Hence, if
$k>2$, the thesis is proved. $\Box$

\section{Appendix} \label{sec:app}

\noindent \emph{Proof of Proposition \ref{prop:disc-stime}.} Part
(i) is an easy consequence of Lemma \ref{lem:disc-limiti}, which
implies that for all $\beta$
$$|\C(n,\alpha,\beta)| \le \sum_{m=0}^{N-1} \ a_1^{(m)} \alpha_m$$
hence the thesis, since $a_1^{(m)} \alpha_m <1$. The same holds
true if any of the terms of the sum is of the type $\bar a_1^{(m)}
\bar \alpha_m$.

\noindent The proof of part (ii) is even more immediate, since
$$\B(n,\alpha,\beta) \le \sum_{m=1}^N \ \bar \beta^m$$
with $\beta^m <1$ for all $m$. Note that having an equality in the
previous relation would mean that each step of the renormalization
procedure is done using $\bar \alpha_m$ and $\bar \beta^m$.

\noindent To study the behaviour of $\res(n,\alpha,\beta)$ let us
start writing $n-1=r_{k_j}+R_1q_1+R_0$ as usual. It is then easy
to realise that $S(n,\alpha,\beta)$ can assume only a finite
number of values. In particular, for each $\beta$, a short
calculation shows that
$$(b_0+1)(1-\beta) \ge S(n,\alpha,\beta) \ge \left\{
\begin{array}{ll}
b_0-a_1 \beta-\beta, & \mbox{ if }\ b_1=0, \\[0.2cm]
b_0-a_1 \beta, & \mbox{ if }\ b_1>0,
\end{array} \right.$$
for all $n \in [r_{k_j}, r_{k_{j+1}}-1]$. A similar result holds
for $S^c(n,\alpha,\beta)$, namely
$$\left\{
\begin{array}{ll}
b_0-a_1 \beta+1-\beta \ge S^c(n,\alpha,\beta) \ge -\beta
(a_1+1-b_0), & \mbox{ if }\ b_1=0, \\[0.2cm]
b_0-a_1 \beta+1 \ge S^c(n,\alpha,\beta) \ge -\beta (a_1-b_0), &
\mbox{ if }\ b_1>0.
\end{array} \right.$$

\noindent Let us first examine $|S(n,\alpha,\beta)|$. We have that
if $b_1=0$ then $(b_0-a_1\beta -\beta)\ge -1$, whereas if $b_1>0$
then $(b_0-a_1\beta)\ge -a_1 \alpha$. Moreover, maximising
$(b_0+1)(1-\beta)$ on $(0,1)$ yields for all $\alpha$ and $\beta$
$$|S(n,\alpha,\beta)| \le \left\{
\begin{array}{ll}
\left( \frac{a_1}{2}+1 \right) \left( 1- \frac{a_1}{2}\ \alpha
\right)\le 1 + \frac{a_1}{4}, & \mbox{ if $a_1$ is even}, \\[0.2cm]
\frac{a_1+1}{2}\ \left( 1- \frac{a_1-1}{2}\ \alpha \right) \le 1 +
\frac{a_1-1}{4}, & \mbox{ if $a_1$ is odd.}
\end{array}
\right.$$ For $|S^c(n,\alpha,\beta)|$, if $b_1=0$ then $(b_0-a_1
\beta+1-\beta) \le 1$, if $b_1>0$ then $(b_0-a_1 \beta+1) \le a_1
\alpha$. To maximise $|S^c(n,\alpha,\beta)|$ we have to consider
separately the cases $b_1=0$ and $b_1>0$.

\noindent If $b_1=0$, then maximising $\beta(a_1+1-b_0)$, for all
$\alpha$ and $\beta$ we get
$$\beta(a_1+1-b_0) \le \left\{
\begin{array}{ll}
\left( \frac{a_1}{2}+1 \right)^2 \ \alpha \le 1+
\frac{a_1+a_2}{4},
& \mbox{ if $a_1$ is even}, \\[0.2cm]
\frac{(a_1+1)(a_1+3)}{4}\  \alpha \le 1+ \frac{a_1+a_2}{4}, &
\mbox{ if $a_1$ is odd.}
\end{array}
\right.$$ But if $b_1=0$, in the next step of the renormalization,
the partial quotient $a_2$ will not appear (see the scheme
(\ref{eq:schema})).

\noindent If $b_1>0$, we have instead to maximise
$\beta(a_1-b_0)$. For all $\alpha$ and $\beta$ it holds
$$\beta(a_1-b_0) \le \left\{
\begin{array}{ll}
\left( \frac{a_1}{2}+1 \right)  \frac{a_1}{2}\ \alpha \le 1+
\frac{a_1}{4}, & \mbox{ if $a_1$ is even}, \\[0.2cm]
\left( \frac{a_1+1}{2}\right)^2   \alpha \le 1+ \frac{a_1}{4}, &
\mbox{ if $a_1$ is odd}.
\end{array}
\right.$$ Applying these inequalities to each renormalization step
yields

$$\res(n,\alpha,\beta) \le \sum_{m=0}^{N} \ \left(
1+\frac{a_1^{(m)}}{4} \right)$$ where we are using $\bar
a_1^{(m)}=a_1^{(m)}-1< a_1^{(m)}$. The thesis of part (iii)
follows using $a_1^{(m)}=a_{m+1}$, where $(a_k)$ are the partial
quotients of $\alpha$. $\Box$

\vskip 0.5cm

\noindent \emph{Proof of Proposition \ref{prop:discr2}.} We need
to show the existence of a $\beta$ such that $d_n(\alpha,\beta)$
has the form we look for. Let $\beta$ satisfy $b_{2k+1}=0$ for all
$k\ge 0$, then following the scheme (\ref{eq:schema}) one works
with couples which are all of the form
$(\alpha_{2m},\beta^{2m-1})$, for $m\ge 0$ (we denote $\beta\equiv
\beta^{-1}$), and there are no inversions. Hence one has
 $\res(n,\alpha,\beta) = \sum_{m=0}^{\frac{ord(n-1)}{2}+1}
S(n(m),\alpha_{2m},\beta^{2m-1})$, where $n(m)$ denotes the number
of iterates at each step, that is $n(0)=n$, $n(1)=j(n-1)+1$, and
in general $n(m)=j(n(m-1)-1)+1$. We now remark that for all
integers of the form $r_{k_j}+R_1 q_1$ there exists $\bar R_0$
such that if $n-1=r_{k_j}+R_1 q_1 +\bar R_0$ then
$S(n(m),\alpha_{2m},\beta^{2m-1})=b_0^{(2m)}(1-\beta^{2m-1})$ for
all $m$, and since this choice depends only on $R_0$, and changing
$R_0$ does not change $j(n-1)$ (see equation (\ref{eq:jr})), we
can choose a sequence $n_k$ such that at each renormalization step
the term $S$ has the chosen value. Moreover, we can show that if
$a_1^{(2m)}\ge 3$ then $b_0^{(2m)}(1-\beta^{2m-1})$ is bigger than
$\frac{a_1^{(2m)}-2}{4}$ for all $m\ge 0$ if
$b_0^{(2m)}=\frac{a_1^{(2m)}}{2}$ or
$b_0^{(2m)}=\frac{a_1^{(2m)}-1}{2}$, according to whether
$a_1^{(2m)}$ is even or odd, respectively. If $a_1^{(2m)}< 3$, we
can choose $\bar R_0$ such that $S\ge 0$. The first part follows.

\noindent For the second part of the proof, we just change a
little bit the argument, choosing $\beta$ such that $b_1>0$ and
$b_{2k}=0$ for all $k\ge 1$. Then by the scheme (\ref{eq:schema}),
we use all the couples $(\alpha_{2m+1},\beta^{2m})$ for all $m\ge
1$, and they are all inverted, since there is only one inversion
at the beginning. Hence $\res(n,\alpha,\beta) = S(n,\alpha,\beta)+
\sum_{m=1}^{\frac{ord(n-1)}{2}+1}
S^c(n(m),\alpha_{2m+1},\beta^{2m})$. As above, we can choose a
subsequence $n_h$ such that for all $m\ge 1$ it holds
$$S^c(n(m),\alpha_{2m+1},\beta^{2m}) =
-\beta^{2m}(a_1^{(2m+1)}-b_0^{(2m+1)}-1)\le
-\frac{a_1^{(2m+1)}-2}{4}$$ by choosing
$b_0^{(2m+1)}=\frac{a_1^{(2m+1)}}{2}$ or
$b_0^{(2m+1)}=\frac{a_1^{(2m+1)}-1}{2}$, according to whether
$a_1^{(2m+1)}$ is even or odd respectively, if $a_1^{(2m)}\ge 3$.
Otherwise, just choose $S^c\le 0$. $\Box$

\vskip 0.5cm

\noindent \emph{Proof of Theorem \ref{teo:disc}.} The upper bound
follows straightforwardly from (\ref{eq:come-pinner}). Moreover,
from Propositions \ref{prop:disc-stime} and \ref{prop:discr2}, we
obtain that if $(\ell_e^2 + \ell_o^2)
>0$ then there exists a subsequence $(n_h)$ such that, denoting
$N_h=ord(n_h)$, for all $\epsilon >0$ there exists $\bar h$ such
that for all $h>\bar h$
$$n_h D^*_{n_h} \ge \sup_\beta |\res(n_h,\alpha,\beta)|-2N_h \ge
\frac{\max \set{\ell_e,\ell_o}-\epsilon}{4} \ \sum_{m=1}^{N_h+1} \
a_m - 3 N_h$$ where we have used the relation
$$\sum_{\stackrel{m=1}{a_{m}\ge 3}}^k\ a_m\ge \sum_{m=1}^k \
a_m - 2k$$ for all $k\ge 1$, which holds also when the sum is
restricted to even or odd indexes. From this it follows
$$\limsup\limits_{n\to \infty}\
\frac{nD_n^*(\alpha)}{\stackrel{ord(n-1)+1}{\sum\limits_{m=1}}
a_m} \ge \frac{\max \set{\ell_e,\ell_o}-\epsilon}{4}$$ for all
$\epsilon >0$. $\Box$

% ----------------------------------------------------------------

\end{document}